\documentclass{amsart}
\title{Selection principles and Baire spaces\footnote{\lowercase{{\bf {\uppercase{K}}ey words and phrases:} {\uppercase{B}}aire space, First category, {\uppercase{B}}anach-{\uppercase{M}}azur game, {\uppercase{M}}enger property, {\uppercase{H}}urewicz property.\\
{{\bf{\uppercase{S}}ubject {\uppercase{C}}lassification:} {\uppercase{P}}rimary 03{\uppercase{E}}99, 54{\uppercase{D}}20, 54{\uppercase{E}}52.}
}}}
\author{by Marion Scheepers}
\newcommand{\sfin}{{\sf S}_{fin}}

\newcommand{\gfin}{{\sf G}_{fin}}
\newcommand{\op}{\mathcal{O}}
\newcommand{\naturals}{{\mathbb N}}
\newtheorem{theorem}{{\bf Theorem }}

\newtheorem{corollary}[theorem]{{\bf Corollary}}

\newcommand{\reals}{{\mathbb R}}
\newcommand{\rationals}{{\mathbb Q}}

\begin{document}
\begin{abstract}
We prove that if $X$ is a separable metric space with the Hurewicz covering property, then the Banach-Mazur game played on $X$ is determined. The implication is not true when ``Hurewicz covering property" is replaced with ``Menger covering property". 
\end{abstract}
\maketitle

\section{Introduction}

The selection principle $\sfin(\mathcal{A},\mathcal{B})$ states that there is for each sequence $(A_n:n\in\naturals)$ with each $A_n\in\mathcal{A}$, a sequence $(B_n:n\in\naturals)$ such that each $B_n \subset A_n$ is finite and $\bigcup_{n\in\naturals}B_n\in\mathcal{B}$. Letting $\op$ denote for the space $X$ the set of all open covers of $X$, the statement $\sfin(\op,\op)$ denotes the Menger property for $X$. Hurewicz \cite{H25}- introduced the Menger property in 1925 and showed that a conjecture of Menger is equivalent to the statement that a metrizable space has the Menger property if, and only if, it is $\sigma$-compact. In 1927 Hurewicz -\cite{H27}-  defined the following stronger version of the Menger property: For each sequence $(\mathcal{U}_n:n\in\naturals)$ of open covers of $X$, there is a sequence $(\mathcal{V}_n:n\in\naturals)$ such that each $x\in X$ is in all but finitely many of the sets $\bigcup\mathcal{V}_n$. This property is said to be the Hurewicz property. In \cite{coc7} it was shown that the Hurewicz property can also be formulated in the form $\sfin(\mathcal{A},\mathcal{B})$, but we will not need that result here. 

It is clear that $\sigma$-compactness implies the Hurewicz property in all finite powers, and that the Hurewicz property implies the Menger property. Fremlin and Miller -\cite{MF}- disproved Menger's Conjecture, thus showing that Menger's property is weaker than $\sigma$-compactness. Numerous examples in the literature show that Menger's property is not necessarily preserved by finite powers. Chaber and Pol -\cite{chaberpol}- showed that the Menger property does not imply the Hurewicz property, and in \cite{coc2} it was shown that the Hurewicz property does not imply $\sigma$-compactness. Also see \cite{TZ}.

This raises the possibility that theorems proven using the hypothesis that some space $X$ is $\sigma$-compact, may be strengthened by proving it using the weaker hypothesis that for all $n$, $X^n$ has Hurewicz's or Menger's property, or that $X$ has Hurewicz's or Menger's property. Several examples of such work can be found in recent literature, for example: \cite{LB}, \cite{ERPol1} and \cite{ERPol2}.
We give such results in this paper in connection with Baire category.

A topological space is said to be Baire if the intersection of any sequence of dense open subsets is a dense set. It is said to be first category if it is a union of countably many nowhere dense sets. If it is not first category, it is said to be second category. In Exercise 25B of \cite{W} the reader is asked to prove the following statement:
\begin{quote}{\tt If $X$ is a $\sigma$-compact space then it is a second category (respectively Baire) space if, and only if $X$ has an element (respectively, dense set of elements) with a compact neighborhood.}
\end{quote}
We examine weakening the hyptohesis ``$X$ is a $\sigma$-compact space".


\section{The Banach-Mazur game and selection principles}

The Banach-Mazur game on $X$, {\sf BM}(X), is played as follows: Players ONE and TWO play an inning per positive integer. In the $n$-th inning ONE chooses a nonempty open set $O_n$; TWO responds with a nonempty open set $T_n\subseteq O_n$. ONE must also obey the rule that for each $n$, $O_{n+1}\subseteq T_n$. A play 
\[
  O_1,\, T_1,\, \cdots,\, O_n,\, T_n,\, \cdots 
\]
is won by TWO if $\bigcap_{n\in\naturals}T_n \neq \emptyset$; otherwise, ONE wins. 

A strategy of a player is a function with domain the set of finite sequences of moves by the opponent, and with values legal moves for the strategy owner. A strategy $\sigma$ for player TWO is said to be a \emph{tactic} if it is of the form $T_n = \sigma(O_n)$ for all $n$. The notion of a tactic for ONE is defined analogously. In \cite{GT} tactics are also called stationary strategies.
The following facts are well-known -\cite{telgarsky}-: 
\begin{enumerate}
\item{$X$ is a Baire space if, and only if, ONE has no winning strategy in {\sf BM}(X).} 
\item{If $X$ is a separable metrizable space such that TWO has a winning strategy in {\sf BM}(X), then X contains a homeomorphic copy of the Cantor set.} 
\item{There are examples of $X$ where neither player has a winning strategy in {\sf BM}(X).} 
\item{If TWO has a winning strategy in {\sf BM}(X), then for each Baire space Y, X$\times$Y is a Baire space.} 
\item{If TWO has a winning strategy in {\sf BM}(X), then all box powers of X are Baire spaces.} 
\end{enumerate} 

Regarding the abovementioned Exercise 25B of \cite{W} one can indeed prove for $\sigma$-compact spaces $X$ that the following statements are equivalent:
\begin{enumerate}
\item{$X$ is a Baire space.}
\item{$X$ has a dense set of points with compact neighborhoods.}
\item{TWO has a winning strategy in {\sf BM}(X).}
\item{TWO has a winning tactic in {\sf BM}(X).}
\end{enumerate}
It follows that in $\sigma$-compact spaces {\sf BM}(X) is determined. We show that this particular consequence of $\sigma$-compactness is not a consequence of the Menger property, but is a consequence of the Hurewicz property.

There is a natural game, $\gfin(\mathcal{A},\mathcal{B})$, that corresponds to the selection principle $\sfin(\mathcal{A},\mathcal{B})$: The game has an inning per positive integer $n$. In the $n$-th inning ONE first chooses an $O_n\in\mathcal{A}$, and TWO then responds with a finite set $T_n\subseteq O_n$. A play $(O_1,\, T_1,\, \cdots,\, O_n,\, T_n,\, \cdots)$ is won by TWO if $\bigcup_{n\in\naturals}T_n\in\mathcal{B}$. Otherwise, ONE wins.

The following equivalence, proved in Theorem 10 of \cite{H25}, is very useful for applications involving the Menger property:
\begin{theorem}[Hurewicz]\label{mengergame} For topological space $X$ the following are equivalent:
\begin{enumerate}
\item{The space has property $\sfin(\op,\op)$.}
\item{ONE has no winning strategy in $\gfin(\op,\op)$.}
\end{enumerate}
\end{theorem}
Below we shall use this equivalence without specifically referencing Theorem \ref{mengergame}.

\begin{theorem}\label{bairemenger} For $X$ be a {\sf T}$_3$-space with $\sfin(\op,\op)$ the following are equivalent:
\begin{enumerate}
\item{TWO has a winning strategy in {\sf BM}(X).}
\item{$D = \{x\in X:x \mbox{ has a neighborhood with compact closure}\}$ is dense in $X$.}
\end{enumerate}
\end{theorem}
{\flushleft{\bf Proof:}} The proof that $(2)\Rightarrow (1)$ does not require that $X$ have property $\sfin(\op,\op)$. We prove $(1)\Rightarrow (2)$ by proving the contrapositive: If $D$ is not dense, then TWO does not have a winning strategy in {\sf BM}(X). Thus: Assume $D$ is not dense, and let F be a strategy for TWO in the game {\sf BM}(X). 

Define a strategy $\sigma$ for ONE of the game $\gfin(\op,\op)$ as follows: First, player ONE of {\sf BM}(X) moves: $B_1$ is a nonempty open set whose closure is disjoint from $D$. TWO's response is $W_1 = F(B_1)$. Now each neighborhood of each x in $W_1$ has a non-compact closure. Choose $x_1\in W_1$. Choose a neighborhood $V_1$ of $x_1$ with $\overline{V_1}\subset W_1$, and an open (in X) cover $\mathcal{A}_1$ of $\overline{V_1}$ such that no finite subset $\mathcal{F}$ of $\mathcal{A}_1$ satisfies $V_1\subset \overline{\bigcup\mathcal{F}}$ (We used $T_3$). Then we define
\[
  \sigma(\emptyset) = \mathcal{A}_1\bigcup\{X\setminus\overline{V_1}\}.
\]
When TWO responds with a finite set $T_1\subset \sigma(\emptyset)$, ONE plans the move $\sigma(T_1)$ as follows: Player ONE of {\sf BM}(X) responds with 
\[
  B_2= W_1\setminus\overline{\bigcup{T_1}},
\]
a nonempty open set. Then TWO of {\sf BM}(X) plays $W_2 =F(B_1,B_2)$. Choose an $x_2\in W_2$ and a neighborhood $V_2$ of $x_2$ with $\overline{V_2}\subset W_2$. Then choose an open (in $X$) cover $\mathcal{A}_2$ of $V_2$ such that no finite subset $\mathcal{F}\subset\mathcal{A}_2$ has $V_2\subset\overline{\bigcup\mathcal{F}}$. Then put
\[
  \sigma(T_1) = \mathcal{A}_2\bigcup\{X\setminus{\overline{V_2}}\}.
\]
When TWO now responds with a finite $T_2\subset\sigma(T_1)$, then ONE plans the move $\sigma(T_1,T_2)$ as follows:  Player ONE of {\sf BM}(X) responds with 
\[
  B_3= W_2\setminus\overline{\bigcup{T_2}},
\]
a nonempty open set. TWO of {\sf BM}(X) applies the strategy $F$ to obtain $W_3 = F(B_1,B_2,B_3)$. Choose an $x_3\in W_3$, and a neighborhood $V_3$ of $x_3$ with $\overline{V_3}\subset W_3$, and then an open (in X) cover $\mathcal{A}_3$ of $\overline{V_3}$ such that for no finite set $\mathcal{F}\subset\mathcal{A}_3$ do we have $\overline{\bigcup\mathcal{F}}\supset V_3$. Then ONE plays
\[
  \sigma(T_1,T_2) = \mathcal{A}_3\bigcup\{X\setminus\overline{V_3}\},
\]
and so on.

Since $X$ has property $\sfin(\op,\op)$, $\sigma$ is not a winning strategy for ONE of $\gfin(\op,\op)$. Thus, consider a $\sigma$-play
\[
  \sigma(\emptyset),\, T_1,\, \sigma(T_1),\, \cdots,\, T_n,\, \sigma(T_1,\cdots,T_n),\, \cdots
\]
lost by ONE. It corresponds to an $F$-play
\[
  B_1,\, F(B_1),\, B_2,\, F(B_1,B_2),\, \cdots,\, B_n,\, F(B_1,\cdots,B_n),\, \cdots
\]
of {\sf BM}(X) where for all $n$ we have $B_{n+1} = F(B_1,\cdots,B_n)\setminus\overline{\bigcup T_n}$. Since ONE lost the $\sigma$-play, the set $\bigcup_{n\in\naturals}T_n$ is an open cover of $X$. For the corresponding play of {\sf BM}(X) we have $\bigcap_{n\in\naturals}B_n \subseteq W_1\setminus \bigcup_{n\in\naturals}(\bigcup T_n) = \emptyset$. Thus, $F$ is not a winning strategy for TWO in {\sf BM}(X). $\diamondsuit$

In general, if TWO has a winning strategy in {\sf BM}(X), then TWO need not have a winning tactic -\cite{debs}. A number of conditions on $X$ that ensures that TWO has a winning strategy if, and only if, TWO has a winning tactic, are known. These include various completeness properties. Theorem \ref{bairemenger} gives another such condition: 
\begin{corollary}\label{twotactic} If $X$ is a {\sf T}$_3$-space with property $\sfin(\op,\op)$, then TWO has a winning strategy in {\sf BM}(X) if, and only if, TWO has a winning tactic.
\end{corollary}

It follows that the {\sf T}$_{3\frac{1}{2}}$-space X of \cite{debs} in which TWO has a winning strategy, but not a winning tactic, in {\sf BM}(X), does not have the Menger property.

\begin{theorem}[CH]\label{mengerinsuff} There is a subspace $X$ of the real line such that:
\begin{enumerate}
\item{$X$ has the property $\sfin(\op,\op)$ in all finite powers, but}
\item{Neither player has a winning strategy in {\sf BM}(X).} 
\end{enumerate}
\end{theorem}
{\flushleft{\bf Proof:}} Consider a Lusin set $X\subset\reals$ which has the property $\sfin(\op,\op)$ in all finite powers. Such is constructed for example in \cite{coc2} or \cite{michael}. We may assume that $X = X+\rationals$. Then for each dense open $E_n\subset X$ there is a dense open $D_n\subset\reals$ with $E_n = X\bigcap D_n$. Since $\reals\setminus D_n$ is nowhere dense, it follows that $X\setminus E_n$ is countable. But then $\bigcap_{n\in\naturals}E_n$ is dense in $X$, showing that $X$ is a Baire space. By the Banach-Oxtoby theorem, ONE has no winning strategy in {\sf BM}(X). Since X contains no subset homeomorphic to the Cantor set, also TWO has no winning strategy in {\sf BM}(X). $\diamondsuit$

We now show that in separable metrizable spaces the Hurewicz property suffices as a replacement for $\sigma$-compactness in the following sense:

\begin{theorem}\label{hurewiczbaire} For $X$ a separable metric space with the Hurewicz property the following are equivalent:
\begin{enumerate}
\item{$X$ is a Baire space.}
\item{$D = \{x\in X:x \mbox{ has a neighborhood with compact closure}\}$ is dense in $X$.}
\item{TWO has a winning strategy in {\sf BM}(X).}
\end{enumerate}
\end{theorem}
{\flushleft{\bf Proof:}} We already have $(2)\Rightarrow (3)$ from Theorem \ref{bairemenger}, and $(3)\Rightarrow(1)$ is folklore. We must prove that $(1)\Rightarrow(2)$. We do this by proving the contrapositive: Assume $D$ is not dense. We will show that ONE has a winning strategy in {\sf BM}(X). The Banach-Oxtoby theorem implies that $X$ is not Baire.

Here is how a winning strategy for ONE is defined. ONE's first move, $\sigma(X)$, is a nonempty open set $O_1\subset X\setminus D$. Since $O_1$ is an {\sf F}$_{\sigma}$ subset of $X$, it has the Hurewicz property also. Fix a metric $d$ on $X$ and choose a countable base $(B_n:n\in\naturals)$ for $O_1$ such that for each $n$, $\overline{B_n}\subset O_1$, $B_n$ has $d$-diameter less than 1, and $\lim_{n\rightarrow\infty}diam(B_n) = 0$ (the latter follows from the Menger property of $O_1$). Now no $\overline{B_n}$ is compact, so we may choose for each $n$ an open (in $O_1$) cover $\mathcal{U}^1_n$ of $\overline{B_n}$ which does not contain any finite set $\mathcal{T}$ with $B_1\subset\overline{\bigcup\mathcal{T}}$. Then for each $n$ the set $\mathcal{U}_n = \{O_1\setminus\overline{B_n}\}\bigcup\mathcal{U}^1_n$ is an open cover of $O_1$. Choose, by the Hurewicz property, for each $n$ a finite set $\mathcal{V}_n\subset\mathcal{U}_n$ such that for each $x\in O_1$, for all but finitely many $n$, $x\in\bigcup\mathcal{V}_n$. We are now ready to define ONE's strategy $\sigma$ further. For each nonempty open set $U\subset O_1$ choose an $n=n(U)$ such that $\overline{B_n}\subset U$, and if $U$ has finite diameter, then $diam_d(B_n)<\frac{1}{2}\cdot diam_d(U)$. When TWO plays an open set $U$, ONE responds with
\[
  \sigma(U) = B_{n(U)}\setminus \overline{\bigcup\mathcal{V}_{n(U)}}.
\]
It is clear that when $U$ is nonempty and open, so is $\sigma(U)$. We must see that $\sigma$ is a winning strategy for ONE. Consider a $\sigma$-play of {\sf BM}(X):
\[
  O_1 = \sigma(X),\, W_1,\, \sigma(W_1),\, W_2,\, \sigma(W_2),\, W_3,\, \cdots
\]
For each $W_k$, put $m_k = n(W_k)$. Then by the definition of ONE's strategy $\sigma(W_1) = B_{m_1}\setminus\overline{\bigcup\mathcal{V}_{m_1}} \supseteq W_2$ and for each $k>1$ $\sigma(W_k) = B_{m_k}\setminus\overline{\bigcup\mathcal{V}_{m_k}}\supseteq W_{k+1}$ and $diam_d(W_{k+1}) < \frac{1}{2}\cdot diam(B_{m_{k-1}})$. This implies that $\{m_k:k\in\naturals\}$ is infinite, so that $\bigcup_{k\in\naturals}\mathcal{V}_{m_k}$ covers $O_1$. It follows that $\bigcap_{k\in\naturals} W_k = \emptyset$, and so ONE wins. $\diamondsuit$

Of course, in Theorem \ref{hurewiczbaire}, we also have the equivalence that TWO has a winning strategy if, and only if, TWO has a winning tactic.

\begin{corollary}\label{bmdetermined} The Banach-Mazur game is determined in separable metric spaces with the Hurewicz property.
\end{corollary}

It is well known that the product of Baire spaces need not be a Baire space again. 
\begin{corollary}\label{baireproducts} If $X$ and $Y$ are Baire spaces and $X$ has the Hurewicz property, then $X\times Y$ is a Baire space.
\end{corollary}

\begin{corollary}\label{bairepowers} If $X$ has the Hurewicz property and is a Baire space, then all powers of $X$ have the Baire property, even in the box topology.
\end{corollary}

However, when $X$ is a separable metric space which has the Baire property and the Hurewicz property, $X^2$ need not have the Hurewicz property. To see this, let C be the Cantor set in $\reals$. Then $Y = \reals\setminus C$ is $\sigma$-compact and Baire. Let $Z\subset C$ be a set with the Hurewicz property in the inherited topology, but  for which $Z\times Z$ does not have the Hurewicz property. The Continuum Hypothesis can be used to find such a subset of the Cantor set - (see the remark following Theorem 2.11 of \cite{coc2}). Put $X = Y\bigcup Z$. Then $X$ is a Baire space and has the Hurewicz property. But the closed subset $Z\times Z$ of $X\times X$ does not have the Hurewicz property, and so $X\times X$ does not have the Hurewicz property.

\section{The game MB(X) and selection principles.}

The game {\sf MB}(X) is played like {\sf BM}(X), except that now ONE wins if $\bigcap_{n\in\naturals}B_n\neq\emptyset$, and TWO wins otherwise.
\begin{theorem}\label{mbandmenger} For $X$ a {\sf T}$_3$-space with $\sfin(\op,\op)$ the following are equivalent:
\begin{enumerate}
\item{ONE has a winning strategy in {\sf MB}(X).}
\item{$D = \{x\in X:x \mbox{ has a neighborhood with compact closure}\}$ is nonempty.}
\end{enumerate}
\end{theorem}
{\flushleft{\bf Proof:}} The proof that $(2)\Rightarrow (1)$ does not require that $X$ has the property $\sfin(\op,\op)$ and uses a standard argument. We prove the contrapositive of $(1)\Rightarrow(2)$: If D is empty then ONE has no winning strategy in {\sf MB}(X). 

Assume $D$ is empty, and let $F$ be a strategy for ONE in the game {\sf MB}(X). Define a strategy $\sigma$ for ONE of the game $\gfin(\op,\op)$ as follows:
Let $B_1=F(X)$ be ONE's move in {\sf MB}(X). Since $\overline{B_1}$ is not compact, choose an open (in $X$) cover $\mathcal{A}_1$ of $\overline{B_1}$ such that no finite $\mathcal{F}\subset\mathcal{A}_1$ has $B_1\subset\overline{\bigcup\mathcal{F}}$, and then define
\[
  \sigma(\emptyset) = \mathcal{A}_1\bigcup\{X\setminus\overline{B_1}\}.
\]
From TWO's response $T_1\subset\sigma(\emptyset)$ in $\gfin(\op,\op)$ define a response $W_1 = B_1\setminus\overline{\bigcup T_1}$ for TWO of {\sf MB}(X), and then apply ONE's strategy $F$ to obtain $B_2 = F(W_1)$. Again, since $\overline{B_2}$ is not compact, choose an open (in $X$) cover $\mathcal{A}_2$ of $\overline{B_2}$ such that no finite $\mathcal{F}\subset\mathcal{A}_2$ has $B_2\subset\overline{\bigcup\mathcal{F}}$, and then define
\[
  \sigma(T_1) = \mathcal{A}_2\bigcup\{X\setminus\overline{B_2}\}.
\]
When TWO responds with the finite set $T_2\subset\sigma(T_1)$, define a response $W_2=B_2\setminus\overline{\bigcup T_2}$ for TWO of {\sf MB}(X). Apply ONE's strategy F in {\sf MB}(X) to obtain $B_3=F(W_1,W_2)$. And since $\overline{B_3}$ is not compact, choose an open (in X) cover $\mathcal{A}_3$ of $\overline{B}_3$ such that there is no finite $\mathcal{F}\subset\mathcal{A}_3$ with $B_3\subset \overline{\bigcup\mathcal{F}}$. Then define
\[
  \sigma(T_1,T_2) = \mathcal{A}_3\bigcup\{X\setminus\overline{B_3}\},
\]
and so on. This defines a strategy for ONE in $\gfin(\op,\op)$. As before, since $X$ has Menger's property $\sfin(\op,\op)$ the strategy $\sigma$ is not a winning strategy for ONE in $\gfin(\op,\op)$. This implies that $F$ is not a winning strategy for ONE in {\sf MB}(X) as follows:

Consider a $\sigma$-play lost by ONE:
\[
  \sigma(\emptyset),\, T_1,\, \sigma(T_1),\, T_2,\, \sigma(T_1,T_2),\, \cdots
\]
The reason ONE lost is that $\bigcup_{n\in\naturals}T_n$ is an open cover of $X$. But the corresponding $F$-play of {\sf MB}(X) is:
\[
  B1,\, W_1 = B_1\setminus\overline{\bigcup T_1},\, B_2,\, W_2 = B_2\setminus\overline{\bigcup T_2},\, \cdots
\]
But then $\bigcap_{n\in\naturals}B_n \subseteq B_1\setminus(\bigcup_{n\in\naturals}\bigcup T_n) = B_1\setminus X = \emptyset$, and so ONE lost.
 $\diamondsuit$

\begin{theorem}[Oxtoby]\label{firstcategory} For a topological space $X$ the following are equivalent:
\begin{enumerate}
\item{TWO has a winning strategy in {\sf MB}(X).}
\item{$X$ is first category in itself.}
\end{enumerate}
\end{theorem}

\begin{theorem}\label{hurewiczcategory} Let $X$ be a separable metric space with the Hurewicz property. Then the following are equivalent:
\begin{enumerate}
\item{$X$ is not first category.}
\item{$D = \{x\in X:x \mbox{ has a neighborhood with compact closure}\}$ is nonempty.}
\item{ONE has a winning strategy in {\sf MB}(X).}
\end{enumerate}
\end{theorem}
{\flushleft{\bf Proof:}} The equivalence of (2) and (3) is in Theorem \ref{mbandmenger}. It is clear from Theorem \ref{firstcategory} that $(3)\Rightarrow (1)$. Thus, we must show $(1)\Rightarrow (2)$. We prove the contrapositive by showing that if $D=\emptyset$, then TWO has a winning strategy in {\sf MB}(X). 
The ideas are as in the proof of Theorem \ref{hurewiczbaire}. Thus, assume that D is empty. Let $(B_n:n\in\naturals)$ enumerate a basis of open sets with finite diameters. By hypothesis no $\overline{B_n}$ is compact. Thus, for each $n$ we may take an open (in X) cover $\mathcal{U}^1_n$ of $\overline{B_n}$ such that no finite subset $\mathcal{F}$ has the property that $B_n\subset\overline{\bigcup\mathcal{F}}$. Then, for each $n$, the set $\mathcal{U}_n = \{X\setminus\overline{B_n}\}\bigcup\mathcal{U}^1_n$ is an open cover of $X$. Since $X$ has the Hurewicz property choose for each $n$ a finite set $\mathcal{V}_n\subset\mathcal{U}_n$ such that for each $x\in X$, for all but finitely many $n$, $x\in\bigcup\mathcal{V}_n$. Now TWO's strategy $\sigma$ is defined as follows: For a nonempty open set $U\subset X$ choose $m=n(U)$ so that $\overline{B_m}\subset U$, and if $U$ has finite diameter, then $diam(B_m)<\frac{1}{2}\cdot diam(U)$. TWO plays 
\[
  \tau(U):= B_{n(U)}\setminus\overline{\bigcup\mathcal{V}_{n(U)}}.
\]
To see that $\tau$ is a winning tactic for TWO, consider a play
\[
  O_1,\, \tau(O_1),\, O_2,\, \tau(O_2),\, \cdots
\]
during which TWO used $\tau$. Then for each $m>1$ we have $O_m\setminus\overline{\bigcup\mathcal{V}_{n(O_m)}} \supset O_{m+1}$ and $diam(O_{m+1})<\frac{1}{2}\cdot diam(O_m)$. It follows that the set $\{n(O_m):m\in\naturals\}$ is infinite, and so $\bigcup_{m\in\naturals}\mathcal{V}_{n(O_m)}$ is a cover of $X$. This implies that $\bigcap_{m\in\naturals}O_m = \emptyset$, and so TWO wins.
$\diamondsuit$

It follows that {\sf MB}(X) is determined in separable metric spaces with the Hurewicz property. It follows that if a subset of the real line has the Hurewicz property but does not contain any perfect set, then it is perfectly meager (since their intersection with any perfect subset of the real line has the Hurewicz property). This gives an alternative proof of Theorem 5.5 of \cite{coc2}.

\bigskip

\begin{center}
  Address
\end{center}

\medskip


\begin{flushleft}
Marion Scheepers                    \\
Department of Mathematics           \\
Boise State University              \\
Boise, Idaho 83725 USA              \\
e-mail: marion@math.boisestate.edu   \\
\end{flushleft}


\begin{thebibliography}{}
\bibitem{LB} L. Babinkostova, \emph{When does the Haver property imply selective screenability?}, {\bf Topology and its Applications} 154 (2007), 1971 - 1979. 
\bibitem{chaberpol} J. Chaber and R. Pol, \emph{A remark on Fremlin–Miller theorem concerning the Menger property and Michael concentrated sets}, preprint of 10.10.2002.
\bibitem{debs} G. Debs, \emph{Strategies gagnantes dans certains jeux topologiques}, {\bf Fundamenta Mathematicae} 126 (1985), 93 - 105
\bibitem{GT} F. Galvin and R. Telg\'arsky, \emph{Stationary strategies in topological games}, {\bf Topology and its Applications} 22 (1986), 51 - 69. 
\bibitem{H25} W. Hurewicz, \emph{\"Uber eine Verallgemeinerung des Borelschen Theorems}, {\bf Mathematische Zeitschrift} 24 (1925), 401 - 425.
\bibitem{H27} W. Hurewicz, \emph{\"Uber Folgen stetiger Funktionen}, {\bf Fundamenta Mathematicae} 9 (1927), 193 - 204.
\bibitem{coc2} W. Just, A.W. Miller, M. Scheepers and P.J. Szeptycki, \emph{Combinatorics of open covers II}, {\bf Topology and its Applications} 73 (1996), 241 - 266.
\bibitem{coc7} Lj.D.R. Ko\v{c}inac and M. Scheepers, \emph{Combinatorics of open covers (VII): Groupability}, {\bf Fundamenta Mathematicae}, 179 (2003), 131 - 155 
\bibitem{michael} E.A. Michael, \emph{Paracompactness and the Lindel\"of property in finite and countable Cartesian products}, {\bf Compositio Mathematicae}  23 (1971), 199 - 244.
\bibitem{MF} A.W. Miller and D.H. Fremlin, \emph{On some properties of Hurewicz, Menger and Rothberger}, {\bf Fundamenta Mathematicae} 129 (1988), 17 - 33.
\bibitem{oxtoby} J.C. Oxtoby, \emph{The Banach-Mazur game and the Banach category theorem}, in: {Contributions to the Theory of Games, Vol. III}, {\bf Annals of Mathematics Studies} 39 (1957), 159 - 163.
\bibitem{ERPol1} E. and R. Pol, \emph{On metric spaces with the Haver property which are Menger spaces}, preprint
\bibitem{ERPol2} E. and R. Pol, \emph{A metric space with the Haver property whose square fails this property}, preprint
\bibitem{telgarsky} R. Telg\'arsky, \emph{Topological games: On the 50th anniversary of the banach-Mazur game}, {\bf The Rocky Mountain Journal of Mathematics} 17:2 (1987), 227 - 276.
\bibitem{TZ} B. Tsaban and L. Zdomsky, \emph{Scales, fields, and a problem of Hurewicz}, {\bf Journal of the European Mathematical Society} (to appear. Preprint locator arXiv:math/0507043v4)
\bibitem{W} S. Willard, \emph{General Topology}, {\bf Addison Wesley Publ. Co.} 1970
\end{thebibliography}
\end{document}